\documentclass{amsart}

\usepackage{amsmath,amsthm,amsfonts,amscd,amssymb}
\usepackage{color}
\usepackage[pdfborder={0 0 0}]{hyperref}
\usepackage{verbatim}

\newtheorem{thm}{Theorem}[section]

\theoremstyle{definition}

\theoremstyle{remark}
\newtheorem{rem}{Remark}[section]

\begin{document}

\title[Effective Kronecker's theorem avoiding algebraic sets]{On an effective variation of Kronecker's approximation theorem avoiding algebraic sets}
\author{Lenny Fukshansky}\thanks{Fukshansky was supported by the NSA grant H98230-1510051 and Simons Foundation grant \#519058.}
\author{Nikolay Moshchevitin}\thanks{Moshchevitin was supported by RNF Grant No. 14-11-00433.}

\address{Department of Mathematics, 850 Columbia Avenue, Claremont McKenna College, Claremont, CA 91711}
\email{lenny@cmc.edu}
\address{Steklov Mathematical Institute of the Russian Academy of Sciences, Gubkina 8, 119991, Moscow, Russia}
\email{moshchevitin@gmail.com}

\subjclass[2010]{11H06, 11G50, 11J68, 11D99}
\keywords{Kronecker's theorem, Diophantine approximation, heights, polynomials, lattices}

\begin{abstract} Let $\Lambda \subset \mathbb R^n$ be an algebraic lattice, coming from a projective module over the ring of integers of a number field $K$. Let $\mathcal Z \subset \mathbb R^n$ be the zero locus of a finite collection of polynomials such that $\Lambda \nsubseteq \mathcal Z$ or a finite union of proper full-rank sublattices of $\Lambda$. Let $K_1$ be the number field generated over $K$ by coordinates of vectors in $\Lambda$, and let $L_1,\dots,L_t$ be linear forms in $n$ variables with algebraic coefficients satisfying an appropriate linear independence condition over $K_1$. For each $\varepsilon > 0$ and $\boldsymbol a \in \mathbb R^n$, we prove the existence of a vector $\boldsymbol x \in \Lambda \setminus \mathcal Z$ of explicitly bounded sup-norm such that
$$\| L_i(\boldsymbol x) - a_i \| < \varepsilon$$
for each $1 \leq i \leq t$, where $\|\ \|$ stands for the distance to the nearest integer. The bound on sup-norm of $\boldsymbol x$ depends on $\varepsilon$, as well as on $\Lambda$, $K$, $\mathcal Z$ and heights of linear forms. This presents a generalization of Kronecker's approximation theorem, establishing an effective result on density of the image of $\Lambda \setminus \mathcal Z$ under the linear forms $L_1,\dots,L_t$ in the $t$-torus~$\mathbb R^t/\mathbb Z^t$.
\end{abstract}

\maketitle

\def\A{{\mathcal A}}
\def\AA{{\mathfrak A}}
\def\B{{\mathcal B}}
\def\C{{\mathcal C}}
\def\D{{\mathcal D}}
\def\E{{\mathcal E}}
\def\F{{\mathcal F}}
\def\x{{\mathcal H}}
\def\I{{\mathcal I}}
\def\J{{\mathcal J}}
\def\K{{\mathcal K}}
\def\kk{{\mathfrak K}}
\def\L{{\mathcal L}}
\def\M{{\mathcal M}}
\def\mm{{\mathfrak m}}
\def\MM{{\mathfrak M}}
\def\O{{\mathcal O}}
\def\OO{{\mathfrak O}}
\def\R{{\mathcal R}}
\def\s{{\mathcal S}}
\def\V{{\mathcal V}}
\def\UU{{\mathfrak U}}
\def\X{{\mathcal X}}
\def\Y{{\mathcal Y}}
\def\Z{{\mathcal Z}}
\def\H{{\mathcal H}}
\def\a{{\mathfrak a}}
\def\b{{\mathfrak b}}
\def\c{{\mathfrak c}}
\def\cee{{\mathbb C}}
\def\pee{{\mathbb P}}
\def\que{{\mathbb Q}}
\def\real{{\mathbb R}}
\def\zed{{\mathbb Z}}
\def\aaa{{\mathbb A}}
\def\ff{{\mathbb F}}
\def\Nn{{\mathbb N}}
\def\kk{{\mathfrak K}}
\def\qbar{{\overline{\mathbb Q}}}
\def\kbar{{\overline{K}}}
\def\taubar{{\overline{\tau}}}
\def\ybar{{\overline{Y}}}
\def\kkbar{{\overline{\mathfrak K}}}
\def\ubar{{\overline{U}}}
\def\eps{{\varepsilon}}
\def\ahat{{\hat \alpha}}
\def\bhat{{\hat \beta}}
\def\gt{{\tilde \gamma}}
\def\h{{\tfrac12}}
\def\ba{{\boldsymbol a}}
\def\be{{\boldsymbol e}}
\def\bei{{\boldsymbol e_i}}
\def\bc{{\boldsymbol c}}
\def\bm{{\boldsymbol m}}
\def\bk{{\boldsymbol k}}
\def\bi{{\boldsymbol i}}
\def\bl{{\boldsymbol l}}
\def\bp{{\boldsymbol p}}
\def\bq{{\boldsymbol q}}
\def\bu{{\boldsymbol u}}
\def\bt{{\boldsymbol t}}
\def\bs{{\boldsymbol s}}
\def\bS{{\boldsymbol S}}
\def\bv{{\boldsymbol v}}
\def\bw{{\boldsymbol w}}
\def\bx{{\boldsymbol x}}
\def\bX{{\boldsymbol X}}
\def\bz{{\boldsymbol z}}
\def\bwy{{\boldsymbol y}}
\def\bY{{\boldsymbol Y}}
\def\bL{{\boldsymbol L}}
\def\baa{{\boldsymbol\alpha}}
\def\bbb{{\boldsymbol\beta}}
\def\bet{{\boldsymbol\eta}}
\def\bxi{{\boldsymbol\xi}}
\def\bo{{\boldsymbol 0}}
\def\bol{{\boldkey 1}_L}
\def\ep{\varepsilon}
\def\p{\boldsymbol\varphi}
\def\q{\boldsymbol\psi}
\def\rank{\operatorname{rank}}
\def\aut{\operatorname{Aut}}
\def\lcm{\operatorname{lcm}}
\def\sgn{\operatorname{sgn}}
\def\spn{\operatorname{span}}
\def\md{\operatorname{mod}}
\def\Norm{\operatorname{Norm}}
\def\dim{\operatorname{dim}}
\def\det{\operatorname{det}}
\def\Vol{\operatorname{Vol}}
\def\rk{\operatorname{rk}}
\def\ord{\operatorname{ord}}
\def\ker{\operatorname{ker}}
\def\div{\operatorname{div}}
\def\Gal{\operatorname{Gal}}
\def\GL{\operatorname{GL}}
\def\id{\operatorname{id}}

\section{Introduction}
\label{intro}

Let $1,\theta_1,\dots,\theta_t$ be $\que$-linearly independent real numbers. The classical approximation theorem of Kronecker then states that the set of points
$$\left\{ (\{ n\theta_1 \}, \dots,\{ n\theta_t \} ) : n \in \zed \right\}$$
is dense in the $t$-torus $\real^t/\zed^t$, where $\{ \cdot \}$ stands for the fractional part of a real number. This result was originally obtained by Kronecker~\cite{kronecker} in 1884, and presents a deep generalization of Dirichlet's 1842 theorem on Diophantine approximation~\cite{dirichlet}; see, for instance,~\cite{hardy} for a detailed exposition of these classical results.

Kronecker's theorem can also be viewed as a statement on density of the image of the integer lattice under collection of linear forms in the torus $\real^t/\zed^t$ (compare to the famous Oppenheim conjecture for quadratic forms). Specifically, if $L_1,\dots,L_t$ are linear forms in $n$ variables with real coefficients $b_{ij}$ so that the set of numbers $1$ and $b_{ij}$ are linearly independent over $\que$, then for any $\eps >0$ and $\ba \in \real^t$ there exists $\bx \in \zed^n$ such that
\begin{equation}
\label{L-kron}
\| L_i(\bx) - a_i \| < \eps\ \forall\ 1 \leq i \leq t,
\end{equation}
where $\|\ \|$ stands for the distance to the nearest integer. A nice survey of a wide variety of results related to Kronecker's theorem is given in~\cite{gonek}. Classical quantitative results in this direction are related to transference theorems for homogeneous and inhomogeneous approximation for the system of linear forms $L_i(\bx)$ (see \cite{jarnik}, Chapter V of \cite{cass:dioph}, \cite{bl}). In particular, these results give effective bounds for the size of the coordinates of the vector $\bx$ in~\eqref{L-kron} under the assumption that there are effective lower bounds for $\max_i \| L_i(\bx) \|$ in the homogeneous case. Some additional effective results can also be found in~\cite{malajovich}, \cite{vorselen}. 
\smallskip

The main goal of this note is somewhat different. We consider linear forms with algebraic coefficients and extend the previously known versions of Kronecker's theorem in three ways:
\begin{enumerate}
\item allow for the approximating vector $\bx$ as in the equation~\eqref{L-kron} above to come from an algebraic lattice $\Lambda$, 
\item exclude vectors from a prescribed union $\Z$ of projective varieties or sublattices not containing this lattice, that is we are interested in approximation vectors $\bx \in \Lambda \setminus \Z$,
\item we obtain effective constants everywhere in our upper bounds.
\end{enumerate}

\noindent
Effective Diophantine avoidance results, exhibiting solutions to a given problem outside of a prescribed algebraic set can be viewed as statements on distribution of such solutions: not only do small solutions exist, they are also sufficiently well distributed so that it is not possible to ``cut them out" by any finite union of varieties. In the recent years, such results were obtained in the general context of Siegel's lemma (also generalizing Faltings' version of Siegel's lemma~\cite{faltings}, \cite{faltings:siegel}, \cite{faltings:siegel_1}) in~\cite{lf:hyper}, \cite{lf:siegel}, \cite{lf:null}, \cite{gaudron}, \cite{gaudron:remond}, \cite{martin}, and in the context of Cassels' theorem on small zeros of quadratic forms and its generalizations in~\cite{lf:smallzeros}, \cite{cfh}, \cite{me_glenn}, \cite{gaudron:remond-1}. We will extend these investigations to Kronecker's theorem. To obtain effective constants in our bounds we use Liouville-type inequalities (see Remark~\ref{schmidt} below for stronger non-effective inequalities of similar type, which can be derived from Schmidt's Subspace Theorem). To give precise statements of our results, we need some notation.
\smallskip

{\it 1. The lattice.} Let $n \geq 1$ be an integer, and for each vector $\bx \in \real^n$ define the sup-norm
$$|\bx| := \max_{1 \leq i \leq n} |x_i|.$$
Let $K$ be a number field of degree $d=r_1+2r_2$ over $\que$, where $r_1$ and $r_2$ are numbers of its real and complex places, respectively, and write $\O_K$ for its ring of integers. Let $1 \leq s \leq w$ be integers, and let $\M \subset K^w$ be an $\O_K$-module such that $\M \otimes_K K \cong K^s$. Write $\D_K(\M)$ for the discriminant of $\M$. Define $\UU_K(\M)$, a fractional $\O_K$-ideal in $K$, to be
\begin{equation}
\label{UKM}
\UU_K(\M) = \left\{ \alpha \in K : \alpha\M \subseteq \O_K^w \right\}.
\end{equation}
We let $\Lambda_K(\M) \subset \real^{wd}$ be the lattice of rank $sd$, which is the image of $\M$ under the standard Minkowski embedding. 
\smallskip

{\it 2. The projective varieties.} Let $m \geq 1$ be an integer. For each $1 \leq i \leq m$, let $\s_i$ be a finite set of homogeneous polynomials in $\real[x_1, \ldots, x_{wd}]$ and $Z(\s_i)$ be its zero set in $\real^{wd}$, that is,
$$Z(\s_i) = \{\bx \in \real^{wd} : P(\bx) = 0 \mbox{ for all } P \in \s_i\}.$$
For the collection $\bS: = \{\s_1, \ldots, \s_m\}$ of finite sets of homogeneous polynomials, define
\begin{equation}
\label{Z_set}
\Z_{\bS} := \bigcup_{i=1}^m Z(\s_i),
\end{equation}
and
\begin{equation}
\label{def_M}
M_{\bS} := \sum_{i=1}^m \max \{\deg P : P \in \s_i\}.
\end{equation}
We allow for the possibility that $\Z_{\bS} = \{ \bo \}$, in which case we take instead $M_{\bS} = 1$. Notice that $\Z_{\bS}$ is an algebraic set, which is a union of a finite collection of projective varieties. Assume that the lattice $\Lambda_K(\M)$ is not contained in the set $\Z_{\bS}$.
\smallskip

{\it 3. The linear forms.} Let $K_1 = K(\Lambda_K(\M))$, i.e. $K_1$ is the number field generated over $K$ by the entries of any basis matrix of the lattice $\Lambda_K(\M)$. Let $B:=(b_{ij})_{1 \leq i \leq t, 1 \leq j \leq wd}$ be a $t \times wd$ matrix with real algebraic entries so that $1,b_{11},\dots,b_{t(wd)}$ are linearly independent over $K_1$, and let $\ell = [E : \que]$ where $E = K_1(b_{11},\dots,b_{t(wd)})$. We will also write $\ell_v = [E_v : \que_v]$ for the local degree of $E$ at every place $v \in M(E)$. Define $t$ linear forms in $wd$ variables
\begin{equation}
\label{L_form}
L_i(x_1,\dots,x_{wd}) = \sum_{j=1}^{wd} b_{ij} x_j \in \real[x_1,\dots,x_{wd}]\ \forall\ 1 \leq i \leq t.
\end{equation}
Our first goal here is to prove the following effective result on density of the image of the set $\Lambda_K(\M) \setminus \Z_{\bS}$ under the linear forms $L_1,\dots,L_t$ in the torus $\real^t / \zed^t$. Let $h$ denote the usual Weil height on algebraic numbers, as well as its extension to vectors with algebraic coordinates; we recall the definition of height along with other necessary notation in Section~\ref{setup}.

\begin{thm} \label{kron1} Let $\ba = (a_1,\dots,a_t) \in \real^t$ and $\eps > 0$. There exist $\bx \in \Lambda_K(\M) \setminus \Z_{\bS}$ and $\bp \in \zed^t$ such that
$$\left| L_i(\bx) - a_i - p_i \right| < \eps$$
and
$$|\bx| \leq \a_K(t,\ell,s)  \left( sd M_{\bS} |\D_K(\M)|^{\frac{s}{2}} \right)^{\kk + 1} \left( (wd)^{\frac{3}{2}} h(B) \right)^{\kk} \c_K(\M,\ell,t)\ \eps^{-\ell+1},$$
where the exponent $\kk = \ell^2(t+1)-\ell$ and the constants are
$$\a_K(t,\ell,s) = 2^{\ell t (\ell-1) + sr_1 \kk + \frac{sd-1}{2}} (t+1)^{3\ell-1} (t!)^{2\ell}$$
and
$$\c_K(\M,\ell,t) = \min \left\{ h(\alpha)^{(\kk+1)sd-1} h(\alpha^{-1})^{\kk} \ :\ \alpha \in \UU_K(\M) \right\}.$$
\end{thm}
\smallskip

\noindent
One special case of Theorem~\ref{kron1} is when~$\Z_{\bS}$ is a union of linear spaces, which means that the point $\bx$ in question is in~$\Lambda_K(\M)$ but outside of a union of sublattices of smaller rank than~$\Lambda_K(\M)$. What if the rank of such sublattices is equal to the rank of $\Lambda_K(\M)$? The next theorem addresses this situation.

\begin{thm} \label{kron2} Let $\ba = (a_1,\dots,a_t) \in \real^t$ and $\eps > 0$. Let $m > 0$ and $\Gamma_1,\dots,\Gamma_m \subset \Lambda_K(\M)$ be proper sublattices of full rank and respective determinants $\D_1,\dots,\D_m$, and let $\D=\D_1 \cdots \D_m$. Then for every $\alpha \in \UU_K(\M)$ there exist $\bx \in \Lambda_K(\M) \setminus \bigcup_{i=1}^m \Gamma_i$ and $\bp \in \zed^t$ such that
$$\left| L_i(\bx) - a_i - p_i \right| < \eps$$
and
$$|\bx| \leq \left( \b_K(t,\ell,s,w) \left( h(\alpha) h(\alpha^{-1}) h(B) \E_{\alpha} \right)^{\kk} \frac{\D\ \eps^{-\ell+1}}{|\D_K(\M)|^{\frac{sm}{2}}} + 1 \right) \E_{\alpha},$$
where the exponent $\kk = \ell^2(t+1)-\ell$, as in Theorem~\ref{kron1}, the constant
$$\b_K(t,\ell,s,w) = 2^{\ell t (\ell-1) + \frac{\kk}{2} + smr_2} (t+1)^{3\ell-1} (t!)^{2\ell} (wd)^{\frac{3\kk}{2}} ,$$
and $\E_{\alpha} =$
\begin{equation}
\label{E-bnd}
\E_{\alpha}(\M,\Gamma_1,\dots,\Gamma_m) := 2^{\frac{sr_1-1}{2}} h(\alpha)^{sd-1} |\D_K(\M)|^{\frac{s}{2}} \left( \sum_{i=1}^m \frac{\D}{\D_i} - m + 1 \right) + \D^{\frac{1}{sd}}.
\end{equation}
\end{thm}

\begin{rem} \label{simpler} The bounds of Theorems~\ref{kron1} and~\ref{kron2} can be recorded in a slightly weaker simplified form as
$$|\bx| \ll \left( \det \Lambda_K(\M) \right)^{\kk + 1} h(B)^{\kk} \c_K(\M,\ell,t)\ \eps^{-\ell+1}$$
and
$$|\bx| \ll \left( \sum_{i=1}^m \frac{\D}{\D_i} \right)^{\kk+2} ( \det \Lambda_K(\M) )^{\kk-m+1} h(B)^{\kk} \c_K(\M,\ell,t)\ \eps^{-\ell+1},$$
respectively, where the constants in the $\ll$ Vinogradov notation depend on the number field $K$ and the integer parameters $t, \ell, s, m$. The expression $\c_K(\M,\ell,t)$ can be viewed as a certain measure of arithmetic complexity of $\M$; in particular, if $\M \subseteq \O_K^w$, then $\c_K(\M,\ell,t) = 1$.
\end{rem}

Here is a sketch of the proofs of Theorems~\ref{kron1} and~\ref{kron2}. We first construct a point~$\bwy \in \Lambda_K(\M)$ of controlled sup-norm, which is outside of~$\Z_{\bS}$ or $\bigcup_{i=1}^m \Gamma_i$, respectively: in the first case, we use the classical Minkowski's Successive Minima Theorem and a version of Alon's Combinatorial Nullstellensatz~\cite{alon} (we use the convenient formulation developed in~\cite{null}), while in the second we employ a recent result of Henk and Thiel~\cite{martin} on points of small norm in a lattice outside of a union of full-rank sublattices. We use $\bwy$ to construct an infinite sequence of points $n \bwy$ satisfying the above conditions, and use an effective version of Kronecker's original theorem to obtain a value of the index $n$ (depending on $\eps > 0$) for which the required inequalities on values of linear forms are satisfied. In other words, our avoidance strategy is to follow the line $n\bwy$ until a necessary point is found. One may wish to use a similar strategy, but following a higher dimensional subspace of the ambient space in the hope of a better bound, however it is difficult to guarantee avoiding our fixed algebraic set with such strategy. A convenient effective version of Kronecker's theorem that we use is worked out in Section~\ref{kron}. It should be remarked that the most important feature of approximation results such as our Theorems~\ref{kron1} and~\ref{kron2} is the exponent on $\eps$ in the bounds for $|\bx|$. As we show, this exponent is the same as in the corresponding bound of the effective version of Kronecker's theorem that we use.

In Section~\ref{setup} we introduce the necessary notation and provide all the details of our setup. We derive an effective version of Kronecker's theorem in Section~\ref{kron}. We then prove Theorem~\ref{kron1} in Section~\ref{proof-kron} and Theorem~\ref{kron2} in Section~\ref{sublattices}.

\bigskip

\section{Notation and setup}
\label{setup}

Let the notation be as in Section~\ref{intro}. Here we introduce some additional notation needed for our algebraic setup. Let the number field $K$ have discriminant $\D_K$, $r_1$ real embeddings $\sigma_1,\dots,\sigma_{r_1}$ of $K$, and $r_2$ conjugate pairs of complex embeddings $\tau_1,\taubar_1,\dots,\tau_{r_2},\taubar_{r_2}$, then $d=r_1+2r_2$. For each $\tau_k$, write $\Re(\tau_k)$ for its real part and $\Im(\tau_k)$ for its imaginary part. Let us write $M(K)$ for the set of all places of $K$, then the archimedean places of $K$ are in correspondence with the embeddings of $K$, and we choose the absolute values $|\ |_{v_1},\dots,|\ |_{v_{r_1+r_2}}$ so that for each $a \in K$
$$|a|_{v_k} = |\sigma_k(a)|\ \forall\ 1 \leq k \leq r_1$$
and
$$|a|_{v_{r_1+k}} = |\tau_k(a)| = \sqrt{ \Re(\tau_k(a))^2 + \Im(\tau_k(a))^2}\ \forall\ 1 \leq k \leq r_2,$$
where $|\ |$ stands for the usual absolute value on $\real$ or $\cee$, respectively. For each $v \in M(K)$, we write $K_v$ for the completion of $K$ at $v$, and for each $n \geq 1$ we define a local norm $|\ |_v : K_v^n \to \real$ by
$$|\ba|_v := \max_{1 \leq j \leq n} |a_j|_v,$$
for each $\ba = (a_1,\dots,a_n) \in K_v^n$. Then the extended Weil height on $K^n$ is given by
$$h(\ba) = \prod_{v \in M(K)} \max \{ 1, |\ba|_v \}^{d_v/d},$$
where $d_v = [K_v : \que_v]$ is the local degree of $K$ at $v$, so that $\sum_{v \mid u} d_v = d$ for each $u \in M(\que)$.

For each integer $n \geq 1$, define the standard Minkowski embedding $\rho^n_K : K^n \to \real^{nd}$ by
$$\rho^n_K(\ba) := \left( \sigma^n_1(\ba),\dots,\sigma^n_{r_1}(\ba),\Re(\tau^n_1(\ba)), \Im(\tau^n_1(\ba)),\dots,\Re(\tau^n_{r_2}(\ba)),\Im(\tau^n_{r_2}(\ba)) \right).$$
We will now use Minkowski embedding to construct lattices from $\O_K$-modules and outline some of their main properties; see~\cite{me_glenn} for further details. Let $1 \leq s \leq w$ be integers, and let $\M \subset K^w$ be an $\O_K$-module such that $\M \otimes_K K \cong K^s$. By the structure theorem for finitely generated projective modules over Dedekind domains (see, for instance~\cite{lang}),
$$\M = \left\{ \sum_{j=1}^s \beta_j \bwy_j : \bwy_j \in \O_K^w,\ \beta_j \in \I_j \right\}$$
for some $\O_K$-fractional ideals $\I_1,\dots,\I_s$ in $K$. By Proposition~13 on p.66 of~\cite{lang}, the discriminant of $\M$ is then
\begin{equation}
\label{module_disc}
\D_K(\M) := \D_K \prod_{j=1}^s \Nn(\I_j)^2,
\end{equation}
where $\Nn(\I_j)$ is the norm of the fractional ideal $\I_j$.

Let $\Lambda_K(\M) := \rho^w_K(\M)$ be an algebraic lattice of rank $sd$ in $\real^{wd}$, then a direct adaptation of Lemma~2 on p.115 of~\cite{lang} implies that the determinant of $\Lambda_K(\M)$ is
\begin{equation}
\label{det_lkm}
\det(\Lambda_K(\M)) = 2^{-sr_2} |\D_K(\M)|^{\frac{s}{2}} = 2^{-sr_2} |\D_K|^{\frac{s}{2}} \prod_{j=1}^s \Nn(\I_j),
\end{equation}
where the last identity follows by~\eqref{module_disc} above. Let $\bx \in \Lambda_K(\M)$, then $\bx = \rho^w_K(\ba)$ for some $\ba \in \M$ and
\begin{equation}
\label{min_bnd}
|\bx| \geq \frac{1}{\sqrt{2}} h(\alpha)^{-1},
\end{equation}
for any $\alpha \in \UU_K(\M)$ by inequality (54) of~\cite{me_glenn}. Let $v \in M(K)$ be an archimedean place, and assume first that it corresponds to a real embedding $\sigma_j$ for some $1 \leq j \leq r_1$, then $|\ba|_v = |\bx|$. On the other hand, if $v$ corresponds to a complex embedding $\tau_j$ for some $1 \leq j \leq r_2$, then $|\ba|_v \leq \left( \sum_{j=1}^{wd} x_j^2 \right)^{1/2} \leq \sqrt{wd}\ |\bx|$. Hence for each $v \mid \infty$,
\begin{equation}
\label{max_bnd}
|\bx| \leq |\ba|_v \leq \sqrt{wd}\ |\bx|.
\end{equation}

Let $L_1,\dots,L_t$ be the linear forms defined in~\eqref{L_form}. For each $1 \leq i \leq t$, we define
$$|L_i|_v = \max_{1 \leq j \leq wd} |b_{ij}|_v,$$
for each place $v \in M(E)$, and define the height of $L_i$ to be
$$h(L_i) = h(b_{i1},\dots,b_{i(wd)}) = \prod_{v \in M(E)} \max \{ 1, |L_i|_v \}^{\ell_v/\ell}.$$
We similarly define the height of the matrix $B$ to be
$$h(B) = h(b_{11},\dots,b_{t(wd)}),$$
then $h(L_i) \leq h(B)$ for all $1 \leq i \leq t$. We are now ready to proceed.
\bigskip

\section{An effective version of Kronecker's theorem}
\label{kron}

In this section we derive an effective version of Kronecker's theorem, which we then use to prove Theorems~\ref{kron1} and~\ref{kron2}. Similar to the setup in the beginning of Section~\ref{intro}, let $1,\theta_1,\dots,\theta_t$ be $\que$-linearly independent real algebraic numbers. For each $1 \leq j \leq t$, let $f_j(x) \in \zed[x]$ be the minimal polynomial of $\theta_j$ of degree $d_j$, $|f_j|$ be the maximum of absolute values of the coefficients of $f_j$, and $A_j$ be the leading coefficient of $f_j$, so $A_j \leq |f_j|$. By Lemma~3.11 of~\cite{waldschmidt},
$$\frac{1}{2^{d_j}}\ |f_j| \leq h(\theta_j)^{d_j} \leq \sqrt{d_j + 1}\ |f_j|,$$
for every $1 \leq j \leq t$. Define $A$ to be the least common multiple of $A_1,\dots,A_t$, so
\begin{equation}
\label{A_bnd}
A \leq \prod_{j=1}^t |f_j| \leq \prod_{j=1}^t (2 h(\theta_j))^{d_j}.
\end{equation}

Let $F=\que(\theta_1,\dots,\theta_t)$ be a number field of degree $e \geq t+1$, then $e \leq \prod_{j=1}^t d_j$.
Let $\theta_{t+1},\dots,\theta_{e-1} \in F$ be such that
$$1 = \theta_0,\theta_1,\dots,\theta_t,\theta_{t+1},\dots,\theta_{e-1}$$
form a $\que$-basis for $F$. Let $\sigma_1,\dots,\sigma_e$ be the embeddings of $F$ into $\cee$. We recall Liouville inequality. For any $\bm = (m_0,\dots,m_t,0,\dots,0) \in \zed^e$,
\begin{equation}
\label{L1}
A^e \prod_{i=1}^e \left| \sum_{j=0}^{e-1} \sigma_i(\theta_j) m_j \right| \geq 1,
\end{equation}
and so
\begin{equation}
\label{L2}
A^e \left( (t+1) \max_{1 \leq i \leq e, 0 \leq j \leq t} |\sigma_i(\theta_j)| \right)^{e-1} |\bm|^{e-1} \| m_1 \theta_1 + \dots + m_t \theta_t \| \geq 1.
\end{equation}
Now observe that
$$\max_{1 \leq i \leq e, 0 \leq j \leq t} |\sigma_i(\theta_j)| \leq \max_{1 \leq j \leq t} h(\theta_j)^{d_j},$$
and so define 
\begin{equation}
\label{const-c1}
\C_1 = \C_1(\theta_1,\dots,\theta_t) := \left( (t+1) \max_{1 \leq j \leq t} h(\theta_j)^{d_j} \right)^{e-1} \prod_{j=1}^t (2 h(\theta_j))^{ed_j}.
\end{equation}
Then for any $\bo \neq \bm \in \zed^t$,
\begin{equation}
\label{mth}
\| m_1 \theta_1 + \dots + m_t \theta_t \| \geq  \C_1^{-1} |\bm|^{-e+1}.
\end{equation}
We will now apply a transference homogeneous-inhomogeneous argument. A transference principle of this sort was first described in Chapter V, \S 4 of~\cite{cass:dioph}; the particular stronger result we are applying here is obtained in~\cite{bl}. Let us write
$$M(\bwy) = \sum_{i=1}^t \theta_i y_i$$
for $\bwy = (y_1,\dots,y_t) \in \zed^t$, and let
$$L_j(x) = \theta_j x,\ 1 \leq j \leq t$$
for $x \in \zed$. Then~\eqref{mth} guarantees that for any $\bo \neq \bwy \in \zed^t$ with $|\bwy| \leq Y$,
$$\| M(\bwy) \| \geq \C_1^{-1} Y^{-(e-1)}.$$
Now applying the transference Lemma 3 of~\cite{bl} to these linear forms, we have that for every $\ba = (a_1,\dots,a_t) \in \real^t$ there exists $x \in \zed$ such that $|x| \leq 2^{-t} ((t+1)!)^2 \C_1 Y^{e-1}$ and
$$\max_{1 \leq j \leq t} \| L_j(x) - a_j \| \leq 2^{-t} ((t+1)!)^2 Y^{-1}.$$
Letting $Q = \left( 2^{t} ((t+1)!)^{-2} Y \right)^{e-1}$, we obtain that
$$\max_{1 \leq j \leq t} \| L_j(x) - \alpha_j \| \leq Q^{-\frac{1}{e-1}}$$
for some $0 \neq x \in \zed$ with $|x| \leq 2^{-et} ((t+1)!)^{2e} \C_1 Q$. Taking $\eps = Q^{-\frac{1}{e-1}}$ immediately yields the following effective version of Kronecker's theorem.

\begin{thm} \label{KR} Let $1,\theta_1,\dots,\theta_t$ be $\que$-linearly independent real algebraic numbers, and let $e = [\que(\theta_1,\dots,\theta_t) : \que]$. Let $\C_1$ be given by~\eqref{const-c1} above, and let $\eps > 0$. Then for any $(a_1,\dots,a_t) \in \real^t$ there exists $q \in \zed \setminus \{0\}$ such that
\begin{equation}
\label{qtheta}
\|q \theta_j - a_j \| \leq \eps,\ 1 \leq j \leq t
\end{equation}
and
$$|q| \leq 2^{-et} ((t+1)!)^{2e} \C_1 \eps^{-e+1}.$$
In particular, if $h(\theta_j) \leq H$ for all $1 \leq j \leq t$ and $\max \{ e, d_1,\dots,d_t \} \leq \ell$, then
$$|q| \leq \left( 2^{\ell t (\ell-1)} (t+1)^{3\ell-1} (t!)^{2\ell} H^{\ell^2(t+1)-\ell} \right) \eps^{-\ell+1}.$$
\end{thm}

\begin{rem} \label{schmidt} Stronger non-effective results can be derived as corollaries of Schmidt's Subspace Theorem. For instance, results discussed in Chapter~6, \S2 of~\cite{schmidt1980} together with the transference principles of Chapter V, \S 4 of~\cite{cass:dioph} and~\cite{bl} imply, for any $\eps > 0$ and $\ba \in \real^t$ under the assumptions of Theorem~\ref{KR}, the existence of $q \in \zed$ satisfying~\ref{qtheta} such that
$$|q| \leq \C'(\delta) \eps^{-t-\delta},$$
for any $\delta > 0$, where the constant $\C'(\delta)$ is non-effective. This would result in the same exponent on $\eps$ in the bounds for $|q|$ in Theorems~\ref{kron1} and~\ref{kron2}, but with non-effective constants.
\end{rem}

\bigskip

\section{Proof of Theorem~\ref{kron1}}
\label{proof-kron}

Here we present the proof of our first result. Since $\Lambda_K(\M) \nsubseteq \Z_{\bS}$, $\Lambda_K(\M) \nsubseteq Z(\s_i)$ for all $1 \leq i \leq m$, and so for each $i$ at least one polynomial $P_i$ in $\s_i$ is not identically zero on $\Lambda_K(\M)$.  Clearly for each $1 \leq i \leq m$, 
$$Z(\s_i) \subseteq Z(P_i) := \left\{ \bx \in \real^{wd} : P_i(\bx) = 0 \right\}.$$
Define
$$P(\bx) = \prod_{i=1}^m P_{i}(\bx),$$
so that $\Z_{\bS} \subseteq Z(P)$ and $\deg(P) \leq M_{\bS}$, while $\Lambda_K(\M) \nsubseteq Z(P)$. Indeed, $Z(P)$ is the union of hypersurfaces $Z(P_1),\dots,Z(P_m)$, and a lattice cannot be covered by a finite union of hypersurfaces unless it is contained in one of them.  We will next construct a point $\bwy \in \Lambda_K(\M)$ of controlled sup-norm such that $P(\bwy) \neq 0$.

Let $V = \spn_{\real} \Lambda_K(\M)$ be the $sd$-dimensional subspace of $\real^{wd}$ spanned by the lattice~$\Lambda_K(\M)$. For a positive real number $\mu$, let us write
$$C_V(\mu) := \left\{ \bx \in V : |\bx| \leq \mu \right\}$$
for the $sd$-dimensional cube with side-length $2\mu$ centered at the origin in $V$, so $C_V(\mu) = \mu C_V(1)$. Let $0 < \lambda_1 \leq \lambda_2 \leq \dots \leq \lambda_{sd}$ be the successive minima of $\Lambda_K(\M)$ with respect to the cube $C_V(1)$. In other words, for each $1 \leq i \leq sd$,
$$\lambda_i := \min \left\{ \mu \in \real_{>0} : \dim_{\real} \spn_{\real} \left( \Lambda_K(\M) \cap C_V(\mu) \right) \geq i \right\}.$$
Let $\bv_1,\dots,\bv_{sd}$ be a collection of linearly independent vectors in $\Lambda_K(\M)$ corresponding to these successive minima, then $|\bv_i| = \lambda_i$. Since the volume of $sd$-dimensional cube $C_V(1)$ is $2^{sd}$, Minkowski's Successive Minima Theorem (see, for instance,~\cite{cass:geom} or~\cite{gruber_lek}) implies that
$$\frac{\det(\Lambda_K(\M))}{(sd)!} \leq \prod_{i=1}^{sd} |\bv_i| \leq \det(\Lambda_K(\M)),$$
where $\frac{1}{\sqrt{2}} h(\alpha)^{-1} \leq |\bv_1| \leq \dots \leq |\bv_{sd}|$, by~\eqref{min_bnd}. This means that
\begin{equation}
\label{vi_bnd}
 |\bv_1| \leq \dots \leq |\bv_{sd}| \leq \left( \sqrt{2} h(\alpha) \right)^{sd-1} \det(\Lambda_K(\M)).
\end{equation}
Let $I(M_{\bS}) = \{ 0,1,2,\dots,M_{\bS} \}$ be the set of the first $M_{\bS}+1$ non-negative integers. For each $\bxi \in I(M_{\bS})^{sd}$, define
$$\bv(\bxi) = \sum_{i=1}^{sd} \xi_i \bv_i,$$
then
\begin{equation}
\label{v_xi_bnd}
|\bv(\bxi)| = \max_{1 \leq j \leq wd} \left| \sum_{i=1}^{sd} \xi_i v_{ij} \right| \leq sd |\bxi| |\bv_{sd}| \leq sd M_{\bS} \left( \sqrt{2} h(\alpha) \right)^{sd-1} \det(\Lambda_K(\M)), 
\end{equation}
by~\eqref{vi_bnd}. Assume that $P(\bv(\bxi)) = 0$ for each $\bxi \in I(M_{\bS})^{sd}$. Then Theorem~4.2 of~\cite{null} implies that $P(\bx)$ must be identically zero on~$V$, which would contradict the fact that $P$ does not vanish identically on $\Lambda_K(\M)$. Hence there must exist some $\bxi \in I(M_{\bS})^{sd}$ such that $P$ does not vanish at the corresponding $\bwy := \bv(\bxi)$, and $|\bwy| \leq sd M_{\bS} \left( \sqrt{2} h(\alpha) \right)^{sd-1} \det(\Lambda_K(\M))$ by~\eqref{v_xi_bnd}. Since $P(\bx)$ is a homogeneous polynomial, it must be true that $P(n \bwy) \neq 0$ for every $n \in \zed_{>0}$. On the other hand, by our construction
$$n \bwy = n \sum_{i=1}^{sd} \xi_i \bv_i \in \spn_{\zed} \left\{ \bv_1,\dots,\bv_{sd} \right\} \subseteq \Lambda_K(\M),$$
and so $\left\{ n \bwy \right\}_{n \in \zed_{>0}}$ gives an infinite sequence of points in $\Lambda_K(\M)$ outside of $\Z_{\bS}$. For each such point, we have
$$L_i(n \bwy) = n L_i(\bwy),\ \forall\ 1 \leq i \leq t.$$
Let us define, for each $1 \leq i \leq t$,
\begin{equation}
\label{theta-i}
\theta_i := L_i(\bwy) = \sum_{j=1}^{wd} b_{ij} y_j \neq 0,
\end{equation}
since $y_j \in K_1$, not all zero, and $b_{ij}$ are $K_1$-linearly independent. Notice that $\theta_1,\dots,\theta_t \in E$, and hence all of them are algebraic numbers of degree $\leq \ell$.

Let $\alpha \in \UU_K(\M)$. Then, by~\eqref{max_bnd}, for each archimedean $v \in M(E)$,
\begin{eqnarray}
\label{th-1}
\max \{1, |\theta_i|_v \} & \leq & \max \{ 1, (wd)^{\frac{3}{2}} |L_i|_v |\bwy| \} \leq (wd)^{\frac{3}{2}} \max \{ 1, |\bwy| \} \max \{ 1, |L_i|_v \} \nonumber \\
& \leq & \sqrt{2}\ (wd)^{\frac{3}{2}} h(\alpha) |\bwy| \max \{ 1, |L_i|_v \},
\end{eqnarray}
by~\eqref{min_bnd}. By~\eqref{v_xi_bnd}, $|\bwy| \leq sd M_{\bS} \left( \sqrt{2} h(\alpha) \right)^{sd-1} \det(\Lambda_K(\M))$, and hence
\begin{equation}
\label{theta_bnd}
\max \{1, |\theta_i|_v \} \leq sd (wd)^{\frac{3}{2}} M_{\bS} \left( \sqrt{2} h(\alpha) \right)^{sd} \det(\Lambda_K(\M)) \max \{ 1, |L_i|_v \}.
\end{equation}
Now suppose $v \in M(E)$ is non-archimedean. Then $\alpha y_j$ is an algebraic integer for each $1 \leq j \leq wd$, and hence $|\alpha y_j|_v = |\alpha|_v |y_j|_v \leq 1$,
meaning that
$$\max \{ 1, |y_1|_v,\dots,|y_{wd}|_v \} \leq \max \{ 1, |\alpha|^{-1}_v \}.$$
Then
\begin{eqnarray}
\label{th-2}
\max \{1, |\theta_i|_v \} & \leq & \max \{1, |L_i|_v \} \max \{ 1, |y_1|_v,\dots,|y_{wd}|_v \} \nonumber \\
& \leq & \max \{ 1, |\alpha^{-1}|_v \} \max \{1, |L_i|_v \},
\end{eqnarray}
for each non-archimedean $v \in M(E)$. Taking a product over all places of $E$, we obtain:
\begin{eqnarray*}
h(\theta_i) & = & \prod_{v \in M(E)} \max \{1, |\theta_i|_v \}^{\frac{\ell_v}{\ell}} = \left( \prod_{v \mid \infty} \max \{1, |\theta_i|_v \}^{\ell_v} \times \prod_{v \nmid \infty} \max \{1, |\theta_i|_v \}^{\ell_v} \right)^{\frac{1}{\ell}} \\
& \leq & sd (wd)^{\frac{3}{2}} M_{\bS} \left( \sqrt{2} h(\alpha) \right)^{sd} \det(\Lambda_K(\M)) h(L_i) \prod_{v \nmid \infty} \max \{ 1, |\alpha^{-1}|_v \}^{\frac{\ell_v}{\ell}} \\
& \leq & sd (wd)^{\frac{3}{2}} M_{\bS} \left( \sqrt{2} h(\alpha) \right)^{sd} h(\alpha^{-1}) \det(\Lambda_K(\M)) h(L_i).
\end{eqnarray*} 
Recalling that $h(L_i) \leq h(B)$ for all $1 \leq i \leq t$, we obtain
\begin{equation}
\label{hti}
h(\theta_i) \leq 2^{\frac{sd}{2}} sd (wd)^{\frac{3}{2}} M_{\bS} h(\alpha)^{sd} h(\alpha^{-1}) \det(\Lambda_K(\M)) h(B),
\end{equation}
for each $1 \leq i \leq t$, where the choice of $\alpha \in \UU_K(\M)$ is arbitrary.

We will now show that $1,\theta_1,\dots,\theta_t$ are $\que$-linearly independent. Suppose not, then there exist $c_0,c_1,\dots,c_t \in \que$, not all zero, such that
$$c_0 = \sum_{i=1}^t c_i \theta_i = \sum_{i=1}^t \sum_{j=1}^{wd} c_i y_j b_{ij},$$
where not all $c_i y_j$ are equal to zero. Recall that $\bwy \in \Lambda_K(\M)$, meaning that coordinates of $\bwy$ are in $K_1$, hence all $c_i y_j$ are in $K_1$. This contradicts the assumption that $1,b_{11},\dots,b_{1(wd)}$ are linearly independent over $K_1$. Hence $1,\theta_1,\dots,\theta_t$ must be linearly independent over~$\que$.

Now let $\ba = (a_1,\dots,a_t) \in \real^t$ and $\eps > 0$, as in the statement of our theorem. Then, by~\eqref{hti} and Theorem~\ref{KR}, there exists $q \in \zed$ and $\bp \in \zed^t$ such that
\begin{eqnarray}
\label{q_bnd}
|q| & \leq & 2^{\ell t (\ell-1)} (t+1)^{3\ell-1} (t!)^{2\ell} \times \nonumber \\
& & \times \left( 2^{\frac{sd}{2}} sd (wd)^{\frac{3}{2}} M_{\bS} h(\alpha)^{sd} h(\alpha^{-1}) \det(\Lambda_K(\M)) h(B) \right)^{\ell^2(t+1)-\ell} \eps^{-\ell+1}
\end{eqnarray}
and
$$\left| q \theta_i - a_i - p_i \right| < \eps\ \forall\ 1 \leq i \leq t.$$
Letting $\bx = q\bwy$, we see that $q\theta_i = L_i(\bx)$ for each $1 \leq i \leq t$ and $|\bx| = |q| |\bwy|$. Combining these observations with~\eqref{v_xi_bnd}, \eqref{q_bnd} and~\eqref{det_lkm} and taking a minimum over all $\alpha \in \UU_K(\M)$ finishes the proof of the theorem.
\bigskip

\section{Proof of Theorem~\ref{kron2}}
\label{sublattices}

Let $\Gamma_1,\dots,\Gamma_m$ be full-rank sublattices of $\Lambda_K(\M)$ of respective determinants $\D_1,\dots,\D_m$. Let $\Omega = \cap_{i=1}^m \Gamma_i$, then $\Omega$ also has full rank and 
$$\D := \D_1 \cdots \D_m \geq \det \Omega.$$
We write $\lambda_i$ for the successive minima of $\Lambda_K(\M)$ and $\lambda_i(\Omega)$ for the successive minima of~$\Omega$. Theorem~1.2 of~\cite{martin} implies that there exists $\bwy \in \Lambda_K(\M) \setminus \bigcup_{i=1}^m \Gamma_i$ such that
$$|\bwy| < \frac{\det \Lambda_K(\M)}{\lambda_1(\Omega)^{sd-1}} \left( \sum_{i=1}^m \frac{\D}{\D_i} - m + 1 \right) + \lambda_1(\Omega).$$
Our first goal is to make this bound more explicit in terms of the parameters of~$\M$. First notice that by Minkowski's Successive Minima Theorem,
$$\lambda_1(\Omega) \leq \left( \prod_{i=1}^{sd} \lambda_i(\Omega) \right)^{1/sd} \leq \left( \det \Omega \right)^{1/sd} \leq \D^{1/sd}.$$
We also need a lower bound on $\lambda_1(\Omega)$. Observe that $\lambda_1(\Omega) \geq \lambda_1$, while $\lambda_1 \geq \frac{1}{\sqrt{2}} h(\alpha)^{-1}$ for any $\alpha \in \UU_K(\M)$, by~\eqref{min_bnd} above. Putting these estimates together, we see that
\begin{equation}
\label{y_O_bnd}
|\bwy| <  \left( \sqrt{2} h(\alpha) \right)^{sd-1} \det \Lambda_K(\M) \left( \sum_{i=1}^m \frac{\D}{\D_i} - m + 1 \right) + \D^{1/sd}
\end{equation}
for any $\alpha \in \UU_K(\M)$.

Since $\bwy \in \Lambda_K(\M)$ and $|\Lambda_K(\M) : \Gamma_i| = \D_i/\det \Lambda_K(\M)$ for each $1 \leq i \leq m$, it follows that
$$\left( g |\Lambda_K(\M) : \Gamma_i| \right) \bwy = \frac{g \D_i}{\det \Lambda_K(\M)} \bwy \in \Gamma_i,$$
for every $g \in \zed$, and hence
$$\left( \frac{g \D_1 \cdots \D_m}{\left( \det \Lambda_K(\M) \right)^m} \right) \bwy = \left( \frac{g \D}{\left( \det \Lambda_K(\M) \right)^m} \right) \bwy \in \Omega,$$
for every $g \in \zed$. Therefore, it must be true that
$$\left( \frac{g \D}{\left( \det \Lambda_K(\M) \right)^m} + 1 \right) \bwy \in \Lambda_K(\M) \setminus \bigcup_{i=1}^m \Gamma_i,$$
for every $g \in \zed$. For brevity, let us write $\D' = \frac{\D}{\left( \det \Lambda_K(\M) \right)^m}$.

From here on, the argument is largely similar to the proof of Theorem~\ref{kron1} above, but with some notable changes. For each $1 \leq i \leq t$, let $\theta_i$ be as in~\eqref{theta-i} for our choice of $\bwy \in \Lambda_K(\M) \setminus \bigcup_{i=1}^m \Gamma_i$ satisfying~\eqref{y_O_bnd} as above, then
$$L_i((g\D'+1) \bwy) = (g\D'+1) \theta_i\ \forall\ 1 \leq i \leq t.$$
Using~\eqref{th-1} with~\eqref{y_O_bnd} instead of~\eqref{v_xi_bnd}, we obtain that $\max \{1, |\theta_i|_v \} \leq$
$$(wd)^{\frac{3}{2}} \left( \left( \sqrt{2} h(\alpha) \right)^{sd} \det \Lambda_K(\M) \left( \sum_{i=1}^m \frac{\D}{\D_i} - m + 1 \right) + \D^{\frac{1}{sd}} \sqrt{2} h(\alpha) \right) \max \{ 1, |L_i|_v \}$$
for all archimedean $v \in M(E)$, while for the non-archimedean $v \in M(E)$,
$$\max \{1, |\theta_i|_v \} \leq \max \{ 1, |\alpha^{-1}|_v \} \max \{1, |L_i|_v \},$$
as in~\eqref{th-2}. Taking the product over all places of $E$, we have for every $1 \leq i \leq t$:
\begin{eqnarray}
\label{th-4}
h(\theta_i) & \leq & (wd)^{\frac{3}{2}} \sqrt{2} h(\alpha) h(\alpha^{-1}) h(B) \times \nonumber \\
&  & \times \left( \left( \sqrt{2} h(\alpha) \right)^{sd-1} \det \Lambda_K(\M) \left( \sum_{i=1}^m \frac{\D}{\D_i} - m + 1 \right) + \D^{\frac{1}{sd}} \right),
\end{eqnarray}
and $1,\theta_1,\dots,\theta_t$ (and hence $1, \D' \theta_1,\dots, \D' \theta_t$) are $\que$-linearly independent by the same reasoning as in the proof of Theorem~\ref{kron1}.

Now let $\ba = (a_1,\dots,a_t) \in \real^t$ and $\eps > 0$, as in the statement of our theorem. Notice that for each $1 \leq i \leq t$,
$$\left| (g\D'+1) \theta_i - a_i - p_i \right| = \left| g (\D' \theta_i) + (\theta_i - a_i) - p_i \right|,$$
for any integers $p_1,\dots,p_t$. Then, applying Theorem~\ref{KR} to approximate the vector $(\theta_1 - a_1,\dots,\theta_t - a_t)$ by the fractional parts of the integer multiples of the vector $(\D' \theta_1,\dots,\D' \theta_t)$, we conclude that there exists $g \in \zed$ and $\bp \in \zed^t$ such that
\begin{eqnarray}
\label{g_bnd}
|g| & \leq & 2^{\ell t (\ell-1)} (t+1)^{3\ell-1} (t!)^{2\ell} \times \nonumber \\
& & \times \left( (wd)^{\frac{3}{2}} \sqrt{2} h(\alpha) h(\alpha^{-1}) h(B) \E_{\alpha}(\M,\Gamma_1,\dots,\Gamma_m) \right)^{\ell^2(t+1)-\ell}\ \eps^{-\ell+1},
\end{eqnarray}
where $\E_{\alpha}(\M,\Gamma_1,\dots,\Gamma_m)$ is as in~\eqref{E-bnd}, and
$$\left| g (\D' \theta_i) + (\theta_i - a_i) - p_i \right| < \eps\ \forall\ 1 \leq i \leq t.$$
Letting $\bx = (g\D'+1)\bwy$, we see that $(g\D'+1) \theta_i = L_i(\bx)$ for each $1 \leq i \leq t$ and $|\bx| = |g\D'+1| |\bwy|$. Combining these observations with~\eqref{y_O_bnd}, \eqref{g_bnd} and~\eqref{det_lkm} finishes the proof of the theorem.

\bigskip

\noindent
{\bf Acknowledgement.} We are grateful for the wonderful hospitality of the Oberwolfach Research Institute for Mathematics: an important part of this work has been done during our Research in Pairs stay at the Institute. We would also like to thank the referee for the helpful remarks.
\bigskip

\bibliographystyle{plain}  
\bibliography{kronecker}        

\end{document}